\documentclass{amsart}

\usepackage{amssymb}

\newtheorem{theorem}{Theorem}[section]

\theoremstyle{definition}

\theoremstyle{remark}

\numberwithin{equation}{section}

\newcommand{\E}{\mathbb{E}}
\newcommand{\R}{\mathbb{R}}
\newcommand{\bx}{\mathbf{x}}
\newcommand{\bt}{\mathbf{t}}
\newcommand{\bu}{\mathbf{u}}

\DeclareMathOperator{\Vol}{Vol}
\DeclareMathOperator{\conv}{conv}
\DeclareMathOperator{\Var}{Var}
\DeclareMathOperator{\Nw}{Nw}
\DeclareMathOperator{\Span}{span}

\newtheorem{corollary}[theorem]{Corollary}
\newtheorem{proposition}[theorem]{Proposition}

\begin{document}

\title[Random Determinants, Mixed Volumes, and Zeros of Gaussian Fields]{Random Determinants, Mixed Volumes of Ellipsoids, and Zeros of Gaussian Random Fields}


\author{Zakhar Kabluchko}
\address{Zakhar Kabluchko, Institute of Stochastics,
Ulm University,
Helmholtzstr.\ 18,
89069 Ulm, Germany}
\email{zakhar.kabluchko@uni-ulm.de}

\author{Dmitry Zaporozhets}
\address{Dmitry Zaporozhets\\
St.\ Petersburg Department of
Steklov Institute of Mathematics,
Fontanka~27,
 191011 St.\ Petersburg,
Russia}
\email{zap1979@gmail.com}
\thanks{The second author is partially supported by RFBR (10-01-00242), NSh-1216.2012.1,  and DFG (436 RUS 113/962/0-1 R) grants.}




\begin{abstract}
Consider a $d\times d$ matrix $M$ whose rows are independent centered non-degenerate Gaussian vectors $\xi_1,\ldots,\xi_d$ with  covariance matrices $\Sigma_1,\dots,\Sigma_d$. Denote by $\mathcal{E}_i$ the location-dispersion ellipsoid of $\xi_i:\mathcal{E}_i=\{\mathbf{x}\in\mathbb{R}^d\,:\, \mathbf{x}^\top\Sigma_i^{-1} \mathbf{x}\leqslant 1\}$. We show that
\begin{equation*}
\mathbb{E}\,|\det M|=\frac{d!}{(2\pi)^{d/2}}V_d(\mathcal{E}_1,\dots,\mathcal{E}_d),
\end{equation*}
where $V_d(\cdot,\dots,\cdot)$ denotes the {\it mixed volume}. We also generalize this result to the case of rectangular matrices. As a direct corollary we get an analytic expression for the mixed volume of $d$ arbitrary ellipsoids in $\mathbb{R}^d$.

As another application, we consider a smooth centered non-degenerate Gaussian random field $X=(X_1,\dots,X_k)^\top:\mathbb{R}^d\to\mathbb{R}^k$. Using Kac-Rice formula, we obtain the geometric interpretation of the intensity of zeros of $X$ in terms of the mixed volume of location-dispersion ellipsoids of the gradients of $X_i/\sqrt{\mathbf{Var} X_i}$. 
This relates zero sets of equations to mixed volumes in a way which is reminiscent of the well-known Bernstein theorem about the number of solutions of the typical system of algebraic equations.
\end{abstract}

\maketitle

\section{Main results}
\subsection{Random determinant and mixed volume of ellipsoids}
Consider independent centered non-degenerate Gaussian random vectors $\xi_1,\dots, \xi_k\in\R^d, k\leqslant d,$ with  covariance matrices $\Sigma_1,\dots,\Sigma_k$. Denote by $\mathcal{E}_i$ the location-dispersion ellipsoid of $\xi_i$:
\begin{equation}\label{1311}
\mathcal{E}_i=\{\bx=(x_1,\dots,x_d)^\top\in\R^d\,:\, \bx^\top\Sigma_i^{-1} \bx\leqslant1\},\quad i=1,\dots,k.
\end{equation}
Denote by $M$ a $k\times d$ matrix whose rows are $\xi_1,\dots,\xi_k$.

\begin{theorem}\label{1645}
It holds that
\begin{equation}\label{1321}
\E\,\sqrt{\det( MM^\top)}=\frac{(d)_k}{(2\pi)^{k/2}\kappa_{d-k}}V_d(\mathcal{E}_1,\dots,\mathcal{E}_k,B,\dots,B),
\end{equation}
where $V_d(\cdot,\dots,\cdot)$ denotes {\it the mixed volume} of $d$ convex bodies in $\R^d$ (see Sect.~\ref{1815} for details),
$B$ is the unit ball in $\R^d$,
$(d)_k=d(d-1)\cdot\dots\cdot (d-k+1)$ is the Pochhammer symbol, and $\kappa_n=\pi^{n/2}/\Gamma(1+n/2)$  denotes the volume of the unit ball in $\R^n$.
\end{theorem}

The left-hand side of \eqref{1321} can be interpreted as the average $k$-dimensional volume of a Gaussian random parallelotope.

\begin{corollary}
In the case $k=d$ it holds that
\begin{equation*}
\E\,|\det M|=\frac{d!}{(2\pi)^{d/2}}V_d(\mathcal{E}_1,\dots,\mathcal{E}_d).
\end{equation*}
\end{corollary}

As another direct corollary we can calculate the mixed volume of $d$ arbitrary ellipsoids in $\R^d$.

\begin{corollary}
If $\mathcal{E}_1,\dots,\mathcal{E}_d$ are arbitrary ellipsoids defined by the symmetric posi\-tive-definite matrices $\Sigma_1,\dots,\Sigma_d$ as in \eqref{1311}, then
\begin{align*}
V_d(\mathcal{E}_1,\dots,\mathcal{E}_d)&=\frac{1}{d!}\prod_{i=1}^d(\det\Sigma_i)^{-1/2}\\
&\times\int_{\R^{d^2}}|\det(x_{ij})|\prod_{i=1}^d\exp\left(-\frac12\bx_i^\top\Sigma_i^{-1}\bx_i\right)dx_{11}\dots dx_{dd},
\end{align*}
where
\begin{equation*}
\bx_i=(x_{i1},\dots,x_{id})^\top.
\end{equation*}

\end{corollary}

The only estimate of the mixed volume of ellipsoids that we know  is due to Barvinok \cite{aB97}. He showed that
\begin{equation*}
\frac{\kappa_d}{3^{(d-1)/2}}\sqrt{D_d(\Sigma_1,\dots,\Sigma_d)}\leqslant V_d(\mathcal{E}_1,\dots,\mathcal{E}_d)\leqslant\kappa_d\sqrt{D_d(\Sigma_1,\dots,\Sigma_d)},
\end{equation*}
where $D_d(\cdot,\dots,\cdot)$ denotes the mixed discriminant of $d$ symmetric $d\times d$ matrices:
\begin{equation*}
D_d(A_1,\dots,A_d)=\frac{1}{d!}\frac{\partial^d}{\partial\lambda_1\dots\partial\lambda_d}\det(\lambda_1A_1+\cdots+\lambda_dA_d)\Big|_{\lambda_1=\cdots=\lambda_d=0}.
\end{equation*}

If $\xi_1,\dots,\xi_k$ are independent standard Gaussian vectors, then $MM^\top$ is a Wishart matrix, and \eqref{1321} turns to (see \cite{sW34},\cite{nG63})
\begin{equation*}
\E\,\sqrt{\det( MM^\top)}=\frac{(d)_k\kappa_d}{(2\pi)^{k/2}\kappa_{d-k}}.
\end{equation*}

\subsection{Intrinsic volumes}
If $\xi_1,\xi_2,\dots, \xi_k\in\R^d, k\leqslant d,$ are identically distributed with the common covariance matrix $\Sigma$ and the  location-dispersion ellipsoid $\mathcal{E}$, then \eqref{1321} turns to
\begin{equation}\label{1330}
\E\,\sqrt{\det( MM^\top)}=\frac{k!}{(2\pi)^{k/2}}V_k(\mathcal{E}),
\end{equation}
where $V_k(\cdot)$ denotes the $k$-th {\it intrinsic volume} of a convex body in $\R^d$:
\begin{equation*}
V_k(K)=\frac{{d\choose k}}{\kappa_{d-k}}V_d(\underbrace{K,\dots,K}_{k\;\rm times},B,\dots,B).
\end{equation*}
The normalization is chosen so that  $V_k(K)$ depends only on $K$ and not on the dimension of the surrounding space, that is, if $\dim K<d$, then the computation of $V_k(K)$ in $\R^d$ leads to the same result as the computation  in the affine span of $K$. In particular, if $\dim K=k$, then $V_k(K)=\Vol_k(K)$, the $k$-dimensional volume of $K$.

It is known that $V_1(K)$ is proportional to the mean width of $K$:
\begin{equation*}
V_1(K)=\frac{d\kappa_d}{2\kappa_{d-1}}w(K).
\end{equation*}
Taking $k=1$ in \eqref{1330}, we obtain that for any centered Gaussian vector $\xi$ with the location-dispersion ellipsoid $\mathcal{E}$ it holds that
\begin{equation}\label{1212}
\E\|\xi\|=\frac{1}{\sqrt{2\pi}}V_1(\mathcal{E}).
\end{equation}
It was pointed out by Mikhail Lifshits that \eqref{1212} is a special case of the following remarkable result of Sudakov.
\subsection{Connection with Sudakov's result}
For our purposes the following finite-dimensional version of Sudakov's theorem suffices. The result in full generality can be found  in~\cite[Proposition 14]{vS76}.
\begin{proposition}
For an arbitrary subset $A\subset\R^d$ we have
\begin{equation}\label{1206}
\E\,\sup_{\bx\in A}\langle \bx,\eta\rangle=\frac1{\sqrt{2\pi}}V_1(\conv(A)),
\end{equation}
where $\eta$ is a standard Gaussian vector in $\R^d$ and $\conv(A)$ is the convex hull of $A$.
\end{proposition}

Let us deduce \eqref{1212} from \eqref{1206}. Consider a matrix $U$ such that $\Sigma=U^{-1}(U^{-1})^\top$ and $U\xi$ is a standard Gaussian vector. Using \eqref{1206} with $A=\mathcal{E}$ and $\eta=U\xi$, we get
\begin{align*}
\E\|\xi\|=\E\sup_{\|\bx\|\leqslant1}\langle \bx,\xi\rangle&=\E\sup_{\|\bx\|\leqslant1}\langle (U^{-1})^\top \bx,U\xi\rangle\\
&=\E\sup_{\|U^\top \bx\|\leqslant1}\langle \bx,U\xi\rangle=\E\sup_{\bx\in\mathcal{E}}\langle \bx,U\xi\rangle=\frac1{\sqrt{2\pi}}V_1(\mathcal E).
\end{align*}

{\bf Open problem:} to obtain a formula generalizing \eqref{1330} and \eqref{1206}.

\subsection{Zeros of Gaussian random fields}
 Let $X(\bt)=(X_1(\bt),\dots,X_k(\bt))^\top:\R^d\to\R^k,\,k\leqslant d,$ be a random field. Following Aza\"is and Wschebor \cite{AW09}, we always assume that the following conditions hold:
\begin{itemize}
\item[(a)] $X$ is Gaussian;
\item[(b)] almost surely, the function $X(\cdot)$ is of class $\mathcal{C}^1$;
\item[(c)] for all $\bt\in\R^d, X(\bt)$ has a non-degenerate distribution;
\item[(d)] almost surely, if $X(\bt)=0$ , then $X'(\bt)$, the Jacobian matrix of $X(\bt)$, has the full rank.
\end{itemize}
Then, almost surely, the level set $X^{-1}(0)$ is a $\mathcal{C}^1$-manifold of dimension $d-k$, and for any Borel set $F$ the Lebesgue measure $\Vol_{d-k}(X^{-1}(0)\cap F)$ is well-defined ($\Vol_0(\cdot)$ denotes the counting measure).

It was shown in \cite[p.~177]{AW09}  that
\begin{equation}\label{1743}
\E \Vol_{d-k}(X^{-1}(0)\cap F)=\int_F\E\left(\sqrt{\det\left(X'(\bt)X'(\bt)^\top\right)}\,\big|\,X(\bt)=0\right)p_{X(\bt)}(0)\,d\bt,
\end{equation}
where $p_{X(\bt)}(\cdot)$ is a density of $X(\bt)$. Thus, the integrand in \eqref{1743} can be interpreted as the intensity of zeros of $X$.

In this paper we consider the special case when $X$ is centered and its coordinates $X_1,\dots,X_k$ are independent. Denote by $\mathcal{E}_i(\bt)$   the location-dispersion ellipsoid of $\nabla [X_i(\bt)/\sqrt{\Var X_i(\bt)}]$.

\begin{theorem}\label{1401}
Let $X$ be a centered random field with independent coordinates defined as above and satisfying conditions (a)-(d). Then
\begin{equation}\label{1037}
\E \Vol_{d-k}(X^{-1}(0)\cap F)=\frac{(d)_k}{(2\pi)^{k}\kappa_{d-k}}\int_FV_d(\mathcal{E}_1(\bt),\dots,\mathcal{E}_k(\bt),B,\dots,B)\,d\bt.
\end{equation}
\end{theorem}

Formula~\eqref{1037} relates zero sets of random equations to mixed volumes. In the case $k=d$ it is therefore reminiscent of the well-known fact from the algebraic geometry which we formulate in the next subsection.

\subsection{Bernstein's theorem}
Consider a complex polynomial in $d$ variables
\begin{equation*}
f(z_1,\dots, z_d)=\sum c_{j_1,\dots,j_d} z_1^{j_1}\dots z_d^{j_d}.
\end{equation*}
The Newton polytope of $f$ is a subset of $\R^d$ defined as
\begin{equation*}
\Nw(f)=\conv\{(j_1,\dots,j_d)\in\mathbb Z^d\,:\,c_{j_1,\dots,j_d}\ne0\}.
\end{equation*}

Let $K_1,\dots, K_d$ be compact convex polytopes in $\R^d$ with vertexes in $\mathbb Z^d$. Consider a system of algebraic equations
\begin{equation*}
\begin{cases}
f_1(z_1,\dots, z_d)=0,\cr
\quad\qquad\dots\cr
f_d(z_1,\dots, z_d)=0,
\end{cases}
\end{equation*}
such that $\Nw(f_i)=K_i$. Bernstein showed \cite{dB75} that for almost all such systems (with respect to the Lebesgue measure in the space of the coefficients of the polynomials) the number of nonzero solutions is equal to
\begin{equation*}
\Vol_0(f^{-1}_1(0)\cap\dots\cap f^{-1}_d(0)\setminus \{\mathbf 0\})=d!V_d(K_1,\dots, K_d).
\end{equation*}

\section{Some essential tools from geometry}\label{1815}
For the basic facts from integral and convex geometry we refer the reader to \cite{BZ80} and \cite{SW08}.

\subsection{Mixed volumes.}
Consider arbitrary convex bodies $K_1,\dots,K_d\subset\R^d$. Min\-kowski showed~\cite{hM11} that $\Vol_d(\lambda_1K_1+\dots+\lambda_dK_d)$, where $\lambda_1,\ldots,\lambda_d\geq 0$, is a homogeneous polynomial of degree $d$ with nonnegative coefficients:
\begin{equation}\label{2138}
\Vol_d(\lambda_1K_1+\dots+\lambda_dK_d)=\sum_{i_1=1}^d\dots\sum_{i_d=1}^d\lambda_{i_1}\dots\lambda_{i_d}V_d(K_{i_1},\dots,K_{i_d}).
\end{equation}
The coefficients $V_d(K_{i_1},\dots,K_{i_d})$ are uniquely determined by the assumption that they are symmetric with respect to permutations of $K_{i_1},\dots,K_{i_d}$. The coefficient $V_d(K_1,\dots,K_d)$ is called the mixed volume of $K_1,\dots,K_d$.
Differentiating \eqref{2138}, we get an alternative definition of the mixed volume:
\begin{equation*}
V_d(K_1,\dots,K_d)=\frac{1}{d!}\frac{\partial^d}{\partial\lambda_1\dots\partial\lambda_d}\Vol_d(\lambda_1K_1+\dots+\lambda_dK_d)\big|_{\lambda_1=\cdots=\lambda_d=0}.
\end{equation*}
For any affine transformation $L$ it holds
\begin{equation}\label{1756}
V_d(LK_{1},\dots,LK_{d})=|\det L| \cdot V_d(K_1,\dots,K_d).
\end{equation}
The following relation can also be stated:
\begin{equation}\label{0220}
\int_{\mathbb{S}^{d-1}}V_{d-1}(P_{\bu}K_1,\dots, P_{\bu}K_{d-1})\,d\bu=\frac{\kappa_{d-1}}{\kappa_d}V_d(K_1,\dots, K_{d-1},B),
\end{equation}
where $d\bu$ is the surface measure on $\mathbb{S}^{d-1}$ normalized to have total mass $1$, and $ P_{\bu}$ denotes the orthogonal projection to the linear hyperplane $\bu^\perp$.

\subsection{Volumes of parallelotopes.} For any $A\subset\R^d$ and $\bx_1,\dots,\bx_k\in\R^d$ denote by $P_{\bx_1,\dots,\bx_k}A$   the orthogonal projection of $A$ to $\Span^\perp\{\bx_1,\dots,\bx_k\}$ (the orthogonal complement of the linear span of $\bx_1,\dots,\bx_k$).
Denote by $H_{\bx_1,\dots,\bx_k}$  the parallelotope generated by the vectors  $\bx_1,\dots,\bx_k$. It is known that
\begin{equation}\label{1738}
\Vol_k(H_{\bx_1,\dots,\bx_k})=\sqrt{\det(AA^\top)},
\end{equation}
where $A$ is a matrix whose rows are $\bx_1,\dots,\bx_k$.

For any $\bx_1,\dots,\bx_d\in\R^d$ and $k=1,\dots, d-1$ it holds that
\begin{equation}\label{0226}
\Vol_d(H_{\bx_1,\dots,\bx_d})=\Vol_k(H_{\bx_1,\dots,\bx_k})\Vol_{d-k}(P_{\bx_{1},\dots,\bx_{k}}H_{\bx_{k+1},\dots,\bx_d}).
\end{equation}

\subsection{Ellipsoids}
There is a bijection $A\mapsto\mathcal E$ between $d\times d$ symmetric positive-definite matrices and $d$-dimensional non-degenerate ellipsoids centered on the origin (see \cite{KVW93} for details):
\begin{equation*}
\mathcal E=\{\bx\in\R^d\,:\, \bx^\top A^{-1} \bx\leqslant1\}.
\end{equation*}
Any non-degenerate linear coordinate transformation of the form $\bx\mapsto L\bx$ is reflected by a change of the corresponding representing matrix $A$ to a matrix $A_L$ given by
\begin{equation}\label{1053}
A_L=LAL^\top.
\end{equation}

Let $\mathcal E'$ be an orthogonal projection of $\mathcal E$ onto an $k$-dimensional subspace with some orthonormal basis $\bx_1,\dots,\bx_k\in\R^d$. Denote by $A'$ a $k\times k$ matrix representing the ellipsoid $\mathcal E'$ in this basis. If $C$ is a $k\times d$ matrix whose rows are $\bx_1,\dots,\bx_k$, then
\begin{equation}\label{1918}
A'=CAC^\top.
\end{equation}

\section{Proofs}
\subsection{Proof of Theorem~\ref{1645}. Case $k=d$}
 We proceed by induction on $d$. First let us assume that $\xi_d$ is a standard Gaussian vector. Denote by $\chi_d$ a random variable having the chi distribution with $d$ degrees of freedom and independent from $\xi_1,\dots,\xi_{d-1}$. Using \eqref{1738} and \eqref{0226} with $k=1$  we get
\begin{align*}
\E\,|\det M|=\E\,\Vol_d(H_{\xi_1,\dots,\xi_d})&=\int_{\mathbb{S}^{d-1}}\E\,\Vol_d(H_{\xi_1,\dots,\xi_{d-1},\chi_d\bu})\,d\bu\\
&=\E\chi_d\int_{\mathbb{S}^{d-1}}\E \Vol_{d-1}(P_{\bu}H_{\xi_1,\dots,\xi_{d-1}})\,d\bu\\
&=\frac{d\kappa_d}{\sqrt{2\pi}\kappa_{d-1}}\int_{\mathbb{S}^{d-1}}\E \Vol_{d-1}(H_{P_{\bu}\xi_1,\dots,P_{\bu}\xi_{d-1}})\,d\bu.
\end{align*}
It follows from \eqref{1918} that $P_{\bu}\xi_i$ has a location-dispersion ellipsoid $P_{\bu}\mathcal E_i$. By the induction assumption,
\begin{equation*}
\E \Vol_{d-1}(H_{P_{\bu}\xi_1,\dots,P_{\bu}\xi_{d-1}})=\frac{(d-1)!}{(2\pi)^{(d-1)/2}}V_{d-1}(P_{\bu}\mathcal E_1,\dots,P_{\bu}\mathcal E_{d-1}).
\end{equation*}
Combining the latter two relations with \eqref{0220}, we obtain
\begin{equation}\label{1611}
\E\,|\det M|=\frac{d!}{(2\pi)^{d/2}}V_{d}(\mathcal E_1,\dots,\mathcal E_{d-1}, B).
\end{equation}

If $\xi_d$ is an arbitrary non-degenerate Gaussian vector, then  there exists a linear transformation $L$ such that $L\xi_d$ is a standard Gaussian vector. It follows from \eqref{1053} that $L\mathcal E_i$ is the location-dispersion ellipsoid of $L\xi_i$, and in particular $L\mathcal E_d=B$.   Applying \eqref{1611} to the matrix $ LM^\top$ and using $\eqref{1756}$, we get
\begin{align*}
\E\,|\det M|=|\det L|^{-1}\,\E\,|\det LM^\top|&=\frac{d!}{(2\pi)^{d/2}}|\det L|^{-1}V_{d}(L\mathcal E_1,\dots,L\mathcal E_{d-1}, B)\\
&=\frac{d!}{(2\pi)^{d/2}}V_{d}(\mathcal E_1,\dots,\mathcal E_{d-1}, \mathcal E_d).
\end{align*}

\subsection{Proof of Theorem~\ref{1645}. Case $k<d$}
Consider a $d\times d$ matrix $M'$ whose first $k$ rows form the matrix $M$ and the last $d-k$ rows are  independent standard Gaussian vectors $\xi_{k+1},\dots,\xi_d$ (independent from $M$). By the previous case,
\begin{equation*}
\E\,|\det M'|=\frac{d!}{(2\pi)^{d/2}}V_{d}(\mathcal E_1,\dots,\mathcal E_{k}, B,\dots, B).
\end{equation*}
On the other hand, by \eqref{0226},
\begin{align*}
\E\,|\det M'|=\E\,\Vol_d(H_{\xi_1,\dots,\xi_d})&=\E\,\Vol_k(H_{\xi_1,\dots,\xi_k})\Vol_{d-k}(P_{\xi_1,\dots,\xi_k}H_{\xi_{k+1},\dots,\xi_d})\\
&=\E\,\sqrt{\det( MM^\top)}\,\E\,\Vol_{d-k}(H_{\eta_{1},\dots,\eta_{d-k}}),
\end{align*}
where $\eta_{1},\dots,\eta_{d-k}$ are independent standard Gaussian vectors in $\R^{d-k}$. By the previous case,
\begin{equation*}
\E\,\Vol_{d-k}(H_{\eta_{1},\dots,\eta_{d-k}})=\frac{(d-k)!}{(2\pi)^{(d-k)/2}}\kappa_{d-k}.
\end{equation*}
Combining the latter three relations completes the proof.

\subsection{Proof of Theorem~\ref{1401}}
First we suppose that $X_j$ has a unit variance: $\Var X_j(\bt)\equiv 1$ for all $j=1,\ldots,k$. Differentiating the relation $\E X_j(\bt)X_j(\bt)=1$ with respect to $t_i$, we obtain
\begin{equation*}
\E\frac{\partial X_j}{\partial t_i}(\bt)X_j(\bt)=0,
\end{equation*}
which, together with the independence of the coordinates of $X$, implies that $X'(\bt)$ and $X(\bt)$ are independent. This means that conditioning on $X(\bt)=0$ in \eqref{1743} may be dropped. To complete the proof of the theorem in the case $\Var X_j(\bt)\equiv 1$ it remains to combine \eqref{1743} with \eqref{1321}.

To cover the general case, it suffices to note that $X_j/\sqrt{\Var X_j}$ has the same zero set as $X_j$.

\bigskip

{\bf Acknowledgments.} We are grateful to Mikhail Lifshits for brining our attention to Sudakov's result.


\providecommand{\bysame}{\leavevmode\hbox to3em{\hrulefill}\thinspace}
\providecommand{\MR}{\relax\ifhmode\unskip\space\fi MR }
\providecommand{\MRhref}[2]{%
  \href{http://www.ams.org/mathscinet-getitem?mr=#1}{#2}
}
\providecommand{\href}[2]{#2}

\end{document}